\documentclass[11pt]{article}
\usepackage{latexsym,amsmath,stackrel,color,amsthm,amssymb,epsfig,graphicx,mathrsfs}
\usepackage{graphicx}
\usepackage{amssymb}
\usepackage[left=1in,top=1in,right=1in,bottom=1in]{geometry}
\usepackage[linktocpage=true]{hyperref}
\usepackage{setspace}
\usepackage{amssymb, amsmath, amsthm, graphicx,mathrsfs}
\usepackage{caption}
\usepackage{titlesec}

\usepackage{comment,caption}
\usepackage{tikz}
\usepackage[compress]{cite}
\makeatletter
\def\thm@space@setup{%
  \thm@preskip=\parskip \thm@postskip=0pt
}
\makeatother

\usepackage{thmtools}

\declaretheoremstyle[%
  spaceabove=6pt,%
  spacebelow=6pt,%
  headfont=\normalfont\itshape,%
  postheadspace=1em,%
  qed=\qedsymbol%
]{mystyle}

\def\qed{\hfill\ifhmode\unskip\nobreak\fi\quad\ifmmode\Box\else\hfill$\Box$\fi}
\def\ite#1{\hfill\break${}$\hbox to 50pt {\quad(#1)\hfill}}
\newtheorem{thm}{Theorem}

\newtheorem{const}{Construction}
\newtheorem{problem}{Problem}

\newtheorem*{remark}{Remark}

\newtheorem{conj}{Conjecture}

\tikzstyle{vertex}=[circle,fill=black,inner sep=2pt]
\tikzstyle{vertrect}=[draw,rectangle,inner sep=2pt]
\tikzstyle{vertdia}=[draw,diamond,inner sep=2pt]

\newcommand{\A}{\mathcal{A}}
\newcommand{\F}{\mathcal{F}}
\newcommand{\B}{\mathcal{B}}

\newcommand{\s}{\mathcal{S}}
\newcommand{\nk}{\binom{[n]}{k}}

\titleformat{\subsection}[runin]
        {\normalfont\bfseries}
        {\thesubsection}
        {0.4em}
        {}
        [.]

\begin{document}
\title{A short note on supersaturation for oddtown and eventown} 
\author{Jason O'Neill \\
Department of Mathematics\\
California State University, Los Angeles\\ 
\tt joneill2@calstatela.edu\footnote{Research partially supported by NSF award DMS-1800332.  \newline\indent
}
} 
\maketitle

\vspace{-10mm}

\begin{abstract}
Given a collection $\A$ of subsets of an $n$ element set, let $\text{op}(\A)$ denote the number of distinct pairs $A,B \in \A$ for which $|A \cap B|$ is odd. For $s \in \{1,2\}$, we prove $\text{op}(\A) \geq s \cdot 2^{\lfloor n/2 \rfloor-1}$ for any collection $\A$ of $2^{\lfloor n/2 \rfloor}+s$ even-sized subsets of an $n$ element set. We also prove $\text{op}(\A) \geq 3$ for any collection $\A$ of $n+1$ odd-sized subsets of an $n$ element set that. Moreover, we show that both of these results are best possible. We then consider larger collections of odd-sized and even-sized sets respectively and explore the connection to minimizing the number of pairwise intersections of size exactly $k-2$ amongst collections of size $k$ subsets from an $n$ element set.     
\end{abstract}

\section{Introduction}

Let $[n]=\{1,2,\ldots,n\}$, $2^{[n]}$ denote the collection of all subsets of $[n]$, and $\nk$ denote the collection of all size $k$ subsets of $[n]$. Given positive integers $k < r$ and a set family $\F \subseteq \binom{[n]}{r}$, let $\partial_k(\F)= \{ A \in \nk: \exists B \in \F, A\subset B\}$ denote the $k$-shadow of the set family $\F$.

Given a family $\A \subseteq 2^{[n]}$, $\A$ follows \textit{oddtown rules} if the sizes of all sets in $\A$ are odd and distinct pairs of sets from $\A$ have even-sized intersections. Berlekamp \cite{BERK} and Graver \cite{GRA} independently proved that the size of a family that satisfies oddtown rules is at most $n$, which is easily seen to be best possible. Alternatively, $\A \subset 2^{[n]}$ of even-sized subsets follows \textit{eventown rules} if all pairs of sets from $\A$ have even-sized intersections. Berlekamp \cite{BERK} and Graver \cite{GRA} independently proved that the size of a family that satisfies eventown rules is at most $2^{\lfloor n/2 \rfloor}$, which is also best possible. For more on these problems, refer to the surveys by Babai and Frankl \cite{BABAIFRANKL} and Frankl and Tokushige \cite{FT}.

The oddtown and eventown problems are foundational results which highlight the linear algebra method \cite{BABAIFRANKL} in extremal combinatorics. There have been numerous extensions and variants of these results in the literature \cite{OV,SV,V1,V2,V0,FO,DFS} and in this paper, we explore their \textit{supersaturation} versions: given slightly more than $n$ odd-sized subsets (resp. $2^{\lfloor n/2 \rfloor}$ even-sized subsets) of an $n$ element set, how many pairs of sets from our collection must have an odd number of elements in common? There has been recent attention in the literature \cite{BL,BW,BDLST,DGS1,DGS2} on supersaturation versions of other foundational problems in extremal set theory. These results collectively show a supersaturation version for the classical Erd{\H o}s-Ko-Rado theorem \cite{EKR} on {\it intersecting families} as well as a supersaturation version for Sperner's theorem \cite{Sperner} on {\it antichains}. In this paper, we provide partial progress towards a supersaturation version to the classical Oddtown and Eventown problems, present several conjectures and discuss the corresponding problem for the Erd{\H o}s-S{\' o}s forbidden intersection problem \cite{ES}.  

For a collection, or set family $\A$, let $\textup{op}(\A)$ denote the number of distinct pairs $A,B \in \A$ for which $|A \cap B|$ is odd. For convenience, suppose $n=2k=4l$ and let $X_1, \ldots, X_l$ be pairwise disjoint sets with $X_i=\{x_{1,i}, x_{2,i}, x_{3,i},x_{4,i}\}$. Given $X_i$, let $A_{2i-1}=\{x_{1,i}, x_{2,i}\}$, $A_{2i}=\{x_{3,i}, x_{4,i}\}$, $B_{2i-1}=\{x_{1,i}, x_{4,i}\}$, and $B_{2i}=\{x_{2,i}, x_{3,i}\}$. Letting $[k] = \{1,2,\ldots, k\}$, define \begin{equation}\label{eq:eventownconst}
\A = \bigg\{ \; \bigcup_{ j \in J} A_j \; : \; J \subseteq [k] \; \bigg\} \hspace{5mm} \text{and} \hspace{5mm} \B = \bigg\{ \; \bigcup_{ j \in J} B_j \; : \; J \subseteq [k] \; \bigg\}.     
\end{equation}
Observe that both $\A$ and $\B$ are extremal (i.e. largest) eventown families and moreover for each $B \in \B \setminus \A$, $\text{op}( \A \cup \{B\} ) = 2^{k-1}$. For $1 \leq s \leq 2^{k}-2^{l}$, let $\A_s$ denote a set family which contains $s$ elements from $\B \setminus \A$ together with $\A$. It then follows that $|\A_s|=2^k+s$ and $\textup{op}(\A_s) = s \cdot 2^{k-1}$. Our main result is that we are able to show this is best possible when $s=1$ and $s=2$:

\begin{thm}\label{thm:maineven}
Let $n\geq 1$ and $s \in \{1,2\}$. If $\A \subset 2^{[n]}$ consist of even-sized subsets with $|\A| \geq 2^{\lfloor n/2 \rfloor}+s$, then $\textup{op}(\A) \geq s \cdot 2^{\lfloor n/2 \rfloor-1}$. 
\end{thm}

Given $3 \leq s \leq 2^{\lfloor n/2 \rfloor}-2^{\lfloor n/4 \rfloor}$, we believe that the set family $\A_s$ minimizes the number of odd-sized pairwise intersections amongst all families of $2^{\lfloor n/2\rfloor}+s$ even-sized set families and conjecture:  

\begin{conj}\label{conj:largerteven}
Let $n \geq 1$ and fix $3 \leq s \leq 2^{\lfloor n/2 \rfloor}-2^{\lfloor n/4 \rfloor}$. If $\A \subset 2^{[n]}$ consist of even-sized subsets with $|\A| \geq 2^{\lfloor n/2 \rfloor}+s$, then $\textup{op}(\A) \geq s \cdot 2^{\lfloor n/2 \rfloor-1}$. 
\end{conj}

The collection of singletons $\F = \{ \{ i \} \}_{ i \in [n] }$ is an extremal oddtown family. Moreover, when $4|n$, another extremal oddtown family is a collection $\F'$ of vertex disjoint $K_4^{(3)}$ - i.e. all triples on four vertices. For $1 \leq s \leq n$, let $\F_s$ be any set family consisting of $\F$ together with exactly $s$ sets from $\F'$. It then follows that $|\F_s| = n+s$ and $\text{op}(\F_s)=3s$. We prove this is best possible for $s=1$:  

\begin{thm}\label{thm:mainodd}
Let $n \geq 1$ and $\A \subset 2^{[n]}$ consist of odd-sized subsets with $|\A| \geq n+1$. Then $\textup{op}(\A) \geq 3$.
\end{thm}

It is also worth noting that there exist extremal examples to Theorem \ref{thm:mainodd} which do not contain an extremal oddtown family. To see this, let $n=5$ and take $\mathcal{X}=\{X_1, \ldots, X_6\}$ where $X_1 = \{1,2,3\}$, $X_2 =\{1,4,5 \}$, $X_3 = \{1,2,4\}$,  $X_4 =\{1,3,5 \}$ , $X_5 = \{1,3,4\}$, and $X_6 =\{1,2,5 \}$. Then $\text{op}(\mathcal{X}) = 3$.
Given $2 \leq s \leq n$, we believe that the set family $\F_s$ minimizes the number of odd-sized pairwise intersections amongst all families of $n+s$ odd-sized set families and conjecture:

\begin{conj}\label{conj:largertodd}
Let $n \geq 1$ and fix $1 \leq s \leq n$. If $\A$ is a collection of odd-sized subsets of an $n$ element set with $|\A| \geq n+s$, then $\textup{op}(\A) \geq 3s$. 
\end{conj}

Conjecture \ref{conj:largertodd} has applications towards a supersaturation version of the forbidden intersection problem. The forbidden intersection problem, posed by Erd{\H o}s and S{\' o}s in 1971 asks for the largest family $\F \subset \nk$ such that $|A \cap B| \neq t$ for all $A,B \in \F$. Frankl and F\"uredi \cite{FF} showed that the largest such family has size $c_{k,t}n^t$ when $2t>k-1$ and size $c_{k,t}n^{k-t-1}$ when $2t \leq k-1$ for some constant $c_{k,t}>0$ depending on $k$ and $t$. In the case when $2t>k-1$, the conjectured extremal family consists of taking all size $k$ subsets contained in a {\it Steiner system}.

For positive integers $t < k<n$, a \textit{Steiner system} $\s = \s(n,k,t) \subset \nk$ is a family where for each $T \in \binom{[n]}{t}$, there exists a unique $F \in \s$ so that $T \subset F$. In the case when $t=k-2$, the conjectured extremal construction is $\partial_k(\s(n,k+1,k-2))$. Provided $n$ satisfies various divisibility conditions, this construction has size $6/(k(k-1)) \cdot \binom{n}{k-2}$.  This follows immediately by noting that given any such family $\F$ and size $(k-3)$ subset $A \subset [n]$, the link $\F(A) = \{ B \in [n] \setminus A: A \cup B \in \F\} $ satisfies oddtown rules and hence $|\F(A)| \leq n-k+3$. Thus, we recover the desired upper bound on $|\F|$ as
\begin{equation}\label{eq:linkoddtown}
\binom{k}{k-3} \cdot |\F| = \sum_{A \in \binom{[n]}{k-3} } |\F(A)| \leq \binom{n}{k-3} \cdot (n-k+3).     
\end{equation}
Given integers $t<k$ and a set family $\F \subseteq \nk$, let $c_{k,t}(\F)$ denote the number of distinct pairs of sets $A,B \in \F$ with $|A \cap B|=t$. Let $\F_0=\partial_k(\s(n,k+1,k-2))$ and consider $A \in \nk$ with $A \notin \F_0$ so that $|A \cap S| = k-1$ for some $S \in \s(n,k+1,k-2)$. Then $c_{k,k-2}(\F_0 \cup \{A\}) = (k-1) + 3\binom{k-1}{2}$. We believe that adding sets of the form  $A \notin \F_0$ so that $|A \cap S| = k-1$ for some $S \in \s(n,k+1,k-2)$ minimizes the number of intersections of size $t=k-2$: 

\begin{problem}\label{problem:tavoiding}
Let $n \geq 1$ and fix $s \geq 1$ and $k \geq 4$. If $\F$ is a collection of size $k$ subsets of an $n$ element set with $|\F| \geq 6/(k(k-1)) \cdot \binom{n}{k-2}+s$, then $c_{k,k-2}(\F) \geq ((k-1) + 3\binom{k-1}{2}) \cdot s$.
\end{problem}

Let $\F$ be a set family as in Problem \ref{problem:tavoiding}. Accepting Conjecture \ref{conj:largertodd} and analyzing \eqref{eq:linkoddtown}, it follows that 
\begin{equation}\label{eq:application}
c_{k,k-2}(\F) \cdot (k-2) \geq  \sum_{A \in \binom{[n]}{k-3}} \text{op}( \F(A) ) \geq s \cdot \binom{k}{3}  \cdot 3.   
\end{equation}

As a result, Conjecture \ref{conj:largertodd} would imply partial progress towards Problem \ref{problem:tavoiding}. In particular, it would imply a gap of a factor less than three in the minimum value of $c_{k,k-2}(\F)$ amongst set families $\F \subseteq \nk$ of size $|\F| \geq 6/(k(k-1)) \cdot \binom{n}{k-2}+s$. The first open case is when $k=4$ and $s=1$ for which a similar argument as in Theorem \ref{thm:mainodd} implies that $c_{4,2}(\F) \geq 6$ whereas Problem \ref{problem:tavoiding} in the case when $k=4$ and $s=1$ asks to show $c_{4,2}(\F) \geq 12$.

\subsection{Organization and Notation}
This paper is organized as follows. In Section \ref{sec:proofofodd}, we prove Theorem \ref{thm:mainodd} and we prove Theorem \ref{thm:maineven} in Section \ref{sec:proofofeven}. In Section \ref{sec:concluding}, we discuss various related problems that may be of independent interest.  For $1 \leq k \leq n$ and $A \subseteq [n]$, let $2^{A}$ denote the collection of subsets of $A$, and $\binom{A}{k}$ denote the size $k$ subsets of $A$. Given non-negative integers $a<b$, let $[a,b] = \{a,a+1, \ldots, b\}$. Moreover, for $A \subseteq [n]$, let $v_A \in \mathbb{F}_2^n$ denote its characteristic vector. For a vector subspace $U \subseteq \mathbb{F}_2^n$, let $U^\perp =\{ v \in \mathbb{F}_2^n : \langle v, u \rangle = 0 \; \forall u \in U  \}$. 

\section{Proof of Theorem \ref{thm:mainodd}} \label{sec:proofofodd}
It suffices to consider $\A = \{A_1, \ldots, A_{n+1}\}$ with $|A_i|$ odd for all $i \in [n+1]$ and to prove that we cannot have $\textup{op}(\A) \leq 2$. Let $v_i := v_{A_i}  \in \mathbb{F}_2^n$ for each $i \in [n+1]$ and consider the vector space
\begin{equation}\label{eq:linearcombo}
V = \{ (\epsilon_1, \ldots, \epsilon_{n+1}) \in \mathbb{F}_2^{n+1}: \epsilon_1 v_1 + \cdots + \epsilon_{n+1} v_{n+1} = 0 \}.  
\end{equation}
Without loss of generality, it suffices to consider three cases: 

\textbf{Case 1:} There exists one such odd intersection pair, which we may assume is $\{A_1,A_2\}$. Let $(\epsilon_1, \ldots, \epsilon_{n+1}) \in V$ be nonzero. Taking inner products with each $v_i$ and $\epsilon_1 v_1 + \cdots + \epsilon_{n+1} v_{n+1}$, we obtain $\epsilon_3 = \epsilon_4 = \cdots = \epsilon_{n+1} = 0$ and $\epsilon_1= \epsilon_2$. If $\epsilon_1= \epsilon_2 = 1$, then $0 = v_1+v_2$ which implies $A_1 = A_2$; a contradiction.  

\textbf{Case 2:} There exist two such odd intersection pairs, which we may assume are $\{A_1,A_2\}$ and $\{A_1,A_3\}$. It then follows that $\{v_2, \ldots, v_{n+1}\}$ is linearly independent. Let $(\epsilon_1, \ldots, \epsilon_{n+1}) \in \mathbb{F}_2^{n+1}$ be a nonzero vector in $V$. If $\epsilon_1 =0$, then $0 = \epsilon_2v_2 + \cdots + \epsilon_{n+1}v_{n+1}$ and hence $\epsilon_2 = \epsilon_3 = \cdots = \epsilon_{n+1} = 0$; a contradiction. Thus, we must have $\epsilon_1=1$. Taking inner products with each $v_i$ and $\epsilon_1 v_1 + \cdots + \epsilon_{n+1} v_{n+1}$, we obtain $\epsilon_1+\epsilon_2 = 0$ and $\epsilon_1+\epsilon_3 = 0$ and $\epsilon_4 = \epsilon_5 = \cdots = \epsilon_{n+1} = 0$. Using that $\epsilon_1=1$, it follows that $\epsilon_1= \epsilon_2 = \epsilon_3 = 1$. As such, $v_1+v_2+v_3 = 0$ which implies that each element in $A_1 \cup A_2 \cup A_3$ is in exactly two of $\{A_1,A_2,A_3\}$. Thus, for each $x \in A_1$, it follows that $x \in A_2$ or $x \in A_3$ but not both and hence we reach a contradiction as
\[ |A_1| = |A_1 \cap A_2| + |A_1 \cap A_3| = 1 + 1 = 0 \pmod{2}. \]

\textbf{Case 3a:} There exist two such odd intersection pairs, which we may assume are $\{A_1,A_2\}$ and $\{A_3,A_4\}$ and $\dim(V)=2$. By the rank-nullity theorem, $ \dim ( \text{span} \{v_1, \ldots, v_{n+1}\}) = n-1$. Moreover, observe that $\{v_1,v_3, v_5,v_6, \ldots, v_{n+1}\}$ is an orthogonal basis. 
Let $(\epsilon_1, \ldots, \epsilon_{n+1}) \in V$ be nonzero. For all $i \geq 5$, taking inner products with $v_i$ and $\epsilon_1 v_1 + \cdots + \epsilon_{n+1} v_{n+1}$, we obtain $\epsilon_i = 0$. For $i \leq 4$, taking inner products with $v_i$ and $\epsilon_1 v_1 + \cdots + \epsilon_{n+1} v_{n+1}$, we recover that $0 = \epsilon_1+\epsilon_2$ and $0= \epsilon_3+\epsilon_4$. Hence, $\epsilon_1=\epsilon_2 \in \{0,1\}$ and $\epsilon_3= \epsilon_4 \in \{0,1\}$. Suppose $\epsilon_1= \epsilon_2 = 1$ and $\epsilon_3=\epsilon_4=0$. Then $0 = \epsilon_1 v_1 + \cdots + \epsilon_{n+1} v_{n+1} = v_1 + v_2 $ and we reach a similar contradiction as in Case 1. Similarly, we cannot have $\epsilon_1=\epsilon_2 = 0$ and $\epsilon_3 = \epsilon_4 =1$. As such, the only possible nonzero vector in $V$ is $(1,1,1,1, 0, \ldots, 0)$ which contradicts the dimension of $V$. 

\textbf{Case 3b:} There exist two such odd intersection pairs, which we may assume is $\{A_1,A_2\}$ and $\{A_3,A_4\}$ and $\dim(V)=1$. The reasoning from Case 3a gives that $V=\{ 0, (1,1,1,1, 0, \ldots, 0)\}$.  By the rank-nullity theorem, $\text{span} \{v_1, \ldots, v_{n+1}\} = \mathbb{F}_2^n$. Clearly any basis for $\mathbb{F}_2^n$ consisting of $n$ vectors from $\{v_1, \ldots, v_{n+1} \}$ cannot contain all of $\{v_1, v_2, v_3, v_4 \}$ as $v_1+ v_2+v_3+v_4 = 0$. Without loss of generality, assume that $\{v_2, \ldots, v_{n+1} \}$ forms a basis for $\mathbb{F}_2^n$ and consider $ U = \text{span}\{ v_3, v_4\}$. It then follows that $\dim(U^\perp) = n-2$ and $\{v_2, v_5, v_6, \ldots, v_{n+1}\}$ forms a basis for $U^\perp$. As such, since $v_1 \in U^\perp$, we have $v_1 = \alpha_2v_2 + \alpha_5v_5 + \cdots + \alpha_{n+1}v_{n+1}$. However, this is a contradiction as $v_1 = v_2+v_3+v_4$ and $\{v_2, \ldots, v_{n+1}\}$ is a basis for $\mathbb{F}_2^n$. 

\section{Proof of Theorem \ref{thm:maineven}}\label{sec:proofofeven}

\subsection{The case of $s=1$}\label{subsec:s1}
Let $\A =\{ A_1, \ldots, A_m \} \subset 2^{[n]}$ be a collection of even-sized subsets where $m=2^{\lfloor n/2 \rfloor}+1$ and suppose $\A' \subset \A$ is a maximal (with respect to $\A$) eventown subfamily of $\A$. If $|\A'| \leq 2^{\lfloor n/2 \rfloor-1}$, then since for each $A \in \A \setminus \A'$, there exists $B \in \A'$ for which $|A \cap B|$ is odd,  \[ \text{op}(\A) \geq |\A \setminus \A'| \geq 2^{\lfloor n/2 \rfloor}+1 -2^{\lfloor n/2 \rfloor-1}> 2^{\lfloor n/2 \rfloor-1}. \] 

If $|\A'| = t > 2^{\lfloor n/2 \rfloor-1}$, then without loss of generality assume that $\A' = \{ A_1, \ldots, A_t\}$. Setting $v_i = v_{A_i} \in \mathbb{F}_2^n$ as the characteristic vector of $A_i$, we consider $W = \text{span} \{ v_1, \ldots, v_t\}$. Then, $W \subseteq W^\perp$ and since $\dim(\mathbb{F}_2^n) = \dim(W) + \dim(W^\perp)$, it follows that $\dim(W) \leq \lfloor n/2 \rfloor$. Further, as $t > 2^{\lfloor n/2 \rfloor-1}$, we get $\dim(W) = \lfloor n/2 \rfloor$. 

For each $A_i$ with $i>t$, we consider the linear functional $\chi_i: W \to \mathbb{F}_2$ given by $w \mapsto \langle w,v_i \rangle$. If $\chi_i(W)=\{1\}$, then $\text{op}(\A' \cup \{A_i\}) =t > 2^{\lfloor n/2 \rfloor-1}$ and we are done. Otherwise, by the first isomorphism theorem, it follows that $W/ \ker(\chi_i) \cong \mathbb{F}_2$ and $\dim( \ker (\chi_i)) = \lfloor n/2 \rfloor-1$. As a result, $|\{ v_1, \ldots, v_t\} \setminus \ker( \chi_i))| \geq t-2^{\lfloor n/2 \rfloor-1}$. Hence, there are at least $t-2^{\lfloor n/2 \rfloor-1}$ sets in $\A'$ which have odd-sized intersection with $A_i$. Since we can do this for all $t < i \leq m$, it follows that 
\begin{equation}\label{eq:manyoddintersections}
\text{op}(\A) \geq \bigg( t-2^{\lfloor n/2 \rfloor-1} \bigg) \cdot (m-t).    
\end{equation}
The result then follows by an optimization problem and considering all $t: 2^{\lfloor n/2 \rfloor-1} <t \leq 2^{\lfloor n/2 \rfloor}$. 

\subsection{The case of $s=2$}
Let $\A =\{A_1, \ldots, A_m\} \subset 2^{[n]}$ be a collection of even-sized subsets with $m = 2^{\lfloor n/2 \rfloor}+2$. Let $\A' \subset \A$ be a maximal (with respect to $\A$) eventown subfamily of $\A$ of maximum size. If $|\A'| \leq 2^{\lfloor n/2 \rfloor-2}+1$, then we may partition $\A$ by $\A = \A_1' \sqcup \A_2' \sqcup \cdots \sqcup \A_l'$ with $l \geq 4$ so that for each $A_i \in \A_i$ and each $j<i$, there exists $B_j \in \A_j$ with $|A_i \cap B_j|$ odd and hence 
\[\text{op}(\A) \geq |\A \setminus \A_1'| + |\A \setminus (\A_1' \sqcup \A_2')| + \cdots + |\A_l'| \geq \frac{3}{4} \cdot 2^{\lfloor n/2 \rfloor} + \frac{1}{2} \cdot 2^{\lfloor n/2 \rfloor} >  2^{\lfloor n/2 \rfloor}. \]

If $|\A'| \in [2^{\lfloor n/2 \rfloor-2}+2, 2^{\lfloor n/2 \rfloor-1}-1]$ or $|\A'| \geq 2^{\lfloor n/2 \rfloor-1}+ 2$, then we are done by \eqref{eq:manyoddintersections}.  

If $|\A'| = 2^{ \lfloor n/2 \rfloor-1}+1 =t$, then without loss of generality assume that $\A' = \{ A_1, \ldots, A_t\}$. Setting $v_i = v_{A_i} \in \mathbb{F}_2^n$ as the characteristic vector of $A_i$, we let $W = \text{span} \{ v_1, \ldots, v_t\}$. As in Section \ref{subsec:s1}, we note that $\dim(W) = \lfloor n/2 \rfloor$ and will consider the linear functional $\chi_i: W \to \mathbb{F}_2$ given by $w \mapsto \langle w,v_i \rangle$  with $W_i:=\ker(\chi_i)$. We then set $\B' = \A \setminus \A'$ and further consider $\B'' \subseteq \B'$ to be the sets $B \in \B'$ for which $\text{op}( \A' \cup \{B\})=1$. We now consider a few cases based on $|\B''|$: 

\textbf{Case 1: $|\B''| \leq 1$ \\ } 
For each $ B \in \B' \setminus B''$, there exist at least two sets in $\A'$ with odd-sized intersections with $B$, so 
\[ \text{op}(\A) \geq |\B''| + 2 \cdot |\B' \setminus \B''| \geq 1 + 2(2^{ \lfloor n/2 \rfloor-1}) >  2^{\lfloor n/2 \rfloor}.  \]

\textbf{Case 2: $2 \leq |\B''| \leq 2^{\lfloor n/2 \rfloor-2}$ \\ } 
As $\B'' \neq \emptyset$, there exists an $i>t$ where $W_i = \ker(\chi_i)$ is so that
$ |\{ v_1, \ldots, v_t\} \setminus W_i|=1.$ 
For each $A_j \in \B' \setminus \B''$, we know that $|\{ v_1, \ldots, v_t\} \setminus W_j| \geq 2$ and in particular $W_i \neq W_j$. As $W_i$ and $W_j$ are vector subspaces of $\mathbb{F}_2^n$ of dimension $\lfloor n/2 \rfloor-1$, it follows that $|W_i \cap W_j| \leq 2^{\lfloor n/2 \rfloor-2}$ and hence $|\{ v_1, \ldots, v_t\} \setminus W_j| \geq |W_i \setminus W_j| \geq 2^{\lfloor n/2 \rfloor-2}$. As such, we get 
\begin{equation}\label{eq:linearalgstability}
\text{op}(\A) \geq |\B''| + (2^{\lfloor n/2 \rfloor-2}) \cdot |\B' \setminus \B''| > 2^{\lfloor n/2 \rfloor}.   
\end{equation}
\textbf{Case 3: $|\B''| \geq 2^{\lfloor n/2 \rfloor-2}$ \\} 
Let $\A'' = \{ A \in \A' : \exists \; B \in \B'' \; \text{so that} \; |A \cap B| \; \text{is odd}\}
.$ We first consider the case where $|\A''| \leq n$ and then will show that we actually necessarily have $|\A''| \leq n$. In this case, there exists $A_i \in \A''$ for which $\text{op}( \B'' \cup \{A_i\}) \geq \frac{1}{n} \cdot 2^{\lfloor n/2 \rfloor-2}$. Let $B_1'', B_2'' \in \B''$ be so that $|A_i \cap B_1''|$ and $|A_i \cap B_2''|$ are both odd. Then $|B_1'' \cap B_2''|$ is odd as otherwise $\A' \setminus \{A_i\} \cup \{B_1'',B_2''\}$ is a eventown subfamily of $\A$ of larger size. As a result, 
\begin{equation}\label{eq:manypairs}
\text{op}(\A) \geq \binom{\frac{1}{n} \cdot 2^{\lfloor n/2 \rfloor-2}}{2} >  2^{\lfloor n/2 \rfloor}.     
\end{equation}
Seeking a contradiction, suppose that $|\A''| \geq n+1$. As each element from $\B''$ has an odd number of elements in common with exactly one set from $\A''$, we may find distinct $X_1, \ldots, X_{n+1} \in \A''$ and $Y_1, \ldots, Y_{n+1} \in \B''$ so that $|X_i \cap Y_i|$ is odd for all $i \in [n+1]$ and $|X_i \cap Y_j|$ is even for all $i \neq j \in [n+1]$. As $\A \subset 2^{[n]}$, this violates the \textit{bipartite oddtown theorem} (see \cite[Exercise 1.8]{BABAIFRANKL}). 

If $|\A'| = 2^{ \lfloor n/2 \rfloor}$ and $\dim(W) = \lfloor n/2 \rfloor$, then we may argue in a similar fashion as above. If $\dim(W) < \lfloor n/2 \rfloor$, then the result follows by \eqref{eq:manyoddintersections}. Thus, we have shown $\text{op}(\A) \geq 2^{\lfloor n/2 \rfloor} $ in all cases.

\section{Concluding Remarks}\label{sec:concluding}

$\bullet$ In this paper, we showed that any collection of $2^{\lfloor n/2 \rfloor}+s$ even-sized subsets of an $n$ element set contain at least $s \cdot 2^{\lfloor n/2 \rfloor-1}$ pairwise intersections of odd size for $s \in \{1,2\}$. It is likely the case that the argument in Section \ref{sec:proofofeven} can be adjusted to handle more small values of $s$ in Conjecture \ref{conj:largerteven}. For instance, the argument in \eqref{eq:manyoddintersections} gives the desired result for all collections $\A$ where $\A$ has a maximal\footnote{with respect to $\A$} eventown subfamily $\A' \subset \A$ with $|\A'| \notin [2^{\lfloor n/2\rfloor-i}, 2^{\lfloor n/2\rfloor-i}+s-1]$ for $1 \leq i < \lfloor n/2 \rfloor$. For other families $\A$, one would like to argue similarly as in \eqref{eq:linearalgstability} and \eqref{eq:manypairs}. It would also be interesting to explore supersaturation for eventown with collections of even-sized subfamilies $\A \subset 2^{[n]}$ with $2\cdot 2^{\lfloor n/2 \rfloor} \lesssim |\A| \lesssim 3\cdot 2^{\lfloor n/2 \rfloor} $, and perhaps taking three extremal eventown families $\A_1$, $\A_2$, and $\A_3$ for which $|\A_i \cap \A_j|$ is as small as possible minimizes $\text{op}(\A)$ in this range.    

Further, the collection $\mathcal{E}$ of {\it all} even sets of $[n]$ has $\textup{op}(\mathcal{E}) \sim 1/2 \cdot \binom{|\mathcal{E}|}{2}$; that is, roughly half of its pairwise intersections have an odd number of elements. One the other hand, $\A \cup \B$ from \eqref{eq:eventownconst} is so that $\textup{op}(\A \cup \B) \sim 1/4 \binom{|\A \cup \B|}{2}$. It would be very interesting to interpolate between these: 

\begin{problem}
Let $\epsilon \in (0,1/2)$ and $n$ be sufficiently large. Determine the maximum value $f(\epsilon)$ so that if $\A \subseteq 2^{[n]}$ is a collection of sets of even size with $|\A| \geq 2^{(1-\epsilon)n}$, then $\textup{op}(\A) \geq f(\epsilon) \binom{|A|}{2}$. 
\end{problem}

$\bullet$ We also considered set families $\A \subset 2^{[n]}$ of odd-sized sets with $|\A| \geq n+1$ and proved that $\text{op}(\A) \geq 3$. It would be interesting to also consider the $k$-uniform version: 

\begin{problem}\label{question:uniformoddtown}
Let $k\geq 3$ be odd. If $\F \subset \nk$ is such that $|\F| \geq n+1$, then is $\text{op}(\F) \geq 4$ for $k=3$ and $\text{op}(\F) \geq 5$ otherwise? 
\end{problem}

The set family $\F_1 \subset \binom{[5]}{3}$ where $\F = \binom{[4]}{3} \cup \{ \{1,3,5\}, \{3,4,5\} \} $  shows that Problem \ref{question:uniformoddtown} is best possible when $k=3$ as $\text{op}(\F_1)=4$. 

For $k \geq 5$ odd, one may take  
$\F_2 = \binom{[k+1]}{k} \cup \binom{[k+2,2k+2]}{k} \cup \{ [k-2] \cup \{k+2,k+3\} \}$ as $\text{op}(\F_2)=5$.

\begin{remark}
{\upshape Since the original arXiV posting of this paper, there has been exciting progress towards Conjecture \ref{conj:largerteven}. Antipov and Cherkashin \cite{AC} used spectral methods to show that $\text{op}(\F) \geq s \cdot 2^{\lfloor n/2 \rfloor-2}$ when $|\F| \geq 2^{n/2}+s$. That is, when $n$ is even, there are at least one-half the number of intersections of odd size as there are in the conjectured extremal family. However, for odd $n$, this is much further from the conjectured extremal family. }   
\end{remark}

\textbf{Acknowledgements:} The author would like to thank Jacques Verstra\"{e}te, Sam Spiro, and Calum Buchanan for engaging discussions regarding this problem. This research was partially supported by NSF award DMS-1800332 while the author was a graduate student.

\end{document}